\documentclass{article}
\usepackage{amsmath, amssymb, amsthm}
\usepackage{amssymb,pb-diagram,pst-all}
\usepackage{amscd}

\usepackage{fancybox}
\usepackage{listings}
\newcommand{\ZZ}{\mathbb Z}
\newcommand{\PP}{\mathbb P}
\newcommand{\QQ}{\mathbb Q}

\newcommand{\mcQ}{\mathcal Q}

\newcommand{\MW}{\mathop {\rm MW}\nolimits}

\newcommand{\mcC}{\mathcal C}

\newtheorem{thm}{Theorem}[section]
\newtheorem{cor}[thm]{Corollary}
\newtheorem{prop}[thm]{Proposition}
\newtheorem{lem}[thm]{Lemma}
\newtheorem{defin}[thm]{Definition}
\newtheorem{exmple}[thm]{Example}
\newtheorem{rem}[thm]{Remark}
\newtheorem{qz}[thm]{Question}

\newcommand{\I}{\mathop {\rm I}\nolimits}

\newenvironment{example}{\begin{exmple}\rm }{\end{exmple}}
\newenvironment{remark}{\begin{rem}\rm }{\end{rem}}

\begin{document}
    
\title{
The matroid structure of  vectors of the Mordell-Weil lattice and the topology of plane quartics and bitangent lines
}
\author{
Ryutaro SATO and Shinzo BANNAI 
}
\maketitle
\section{Introduction}
In this paper, we study the embedded topology of plane curves. We are interested in the following situation. Let $\mcC_1, \mcC_2\subset \PP^2$ be plane curves. Then $(\PP^2, \mcC_1)$ and $(\PP^2, \mcC_2)$ form a Zariski-pair if the following conditions are satisfied
\begin{enumerate}
\item There exist tubular neighborhoods $T(\mcC_i)$ of $\mcC_i$ ($i=1,2$) such that the pairs $(T(\mcC_1),\mcC_1)$ and $(T(\mcC_2),\mcC_2)$ are homeomorphic as pairs. 
\item The pairs $(\PP^2, \mcC_1)$ and $(\PP^2, \mcC_2)$ are not homeomorphic as pairs.
\end{enumerate} 
The notion of a Zariski-pair was first defined in \cite{artal94} by E. Artal--Bartolo and has been an object of interest to many mathematicians.
The key in studying Zariski pairs is finding a suitable method to distinguish the curves. Many invariants have been used, such as  the fundamental groups of the complements $\pi_1(\PP^2\setminus \mcC_i)$, the Alexander polynomials $\Delta_{\mcC_i}(t)$ and the existence/non-existence of certain Galois covers branched along $\mcC_i$ (see \cite{act} for a survey on these topics). More recently, newer types of invariants such as \lq\lq linking invariants" and \lq\lq splitting invariants" have been developed in studying reducible plane curves (\cite{bannai16,  benoit-jb, shirane16}). 
However, as the number of irreducible components of $C_i$ increases,  these invariants become more increasingly complex, and it becomes hard to grasp the situation clearly. Hence, we are especially interested in formulating a method in order to present the differences in the curves and the classification comprehensively. 

An attempt at this was done in \cite{bannai-tokunaga15},\cite{bbst} where the second author together with colleagues considered invariants of subsets of the set of irreducible components. This approach proved to be effective and was able to produce new examples of Zariski pairs. However the examples produced were relatively simple, maybe too simple,  to appreciate the usefulness of the approach fully. In this paper, we introduce the terminology of {\it matroids} into our setting in order to make the results more accessible to a wider audience and also to present  more complex examples to demonstrate the usefulness of considering subarrangements more fully. 

We introduce some notation to explain the kind of arrangements that we will study.  Let $Q$ be a smooth quartic curve and $z_o\in Q$ be a general point of $Q$. It is known that a rational elliptic surface $S_{Q, z_0}$ can be associated to $Q$ and $z_o$ as follows (see \cite{tokunaga14, bannai-tokunaga15} for details): Let $\tilde{f}_Q:\tilde{S}_{Q}\rightarrow \PP^2$ be the double cover of $\PP^2$ branched along $Q$, and let $\mu: S_Q\rightarrow \tilde{S}_Q$ be the canonical resolution of singularities. Also, let  $\Lambda_{z_o}$ be the pencil of lines through $z_o$. Then the inverse image $\overline{\Lambda}_{z_o}$ of $\Lambda_{z_o}$  in $\overline{S}_Q$ gives rise to a pencil of curves with genus 1. Next,  the base points of $\overline{\Lambda}_{z_o}$ can be resolved by two consecutive blow-ups, whose composition is denoted by  $\nu_{z_o}:S_{Q, z_o}\rightarrow \overline{S}_Q$. The morphism $\phi_{z_o}: S_{Q, z_o}\rightarrow \PP^1$ induced by $\overline{\Lambda}_{z_o}$ gives a genus 1 fibration, and the exceptional divisor of the second blow-up in $\mu_{z_o}$ gives a section denote by $O$. Hence, we have an elliptic surface $ \phi_{z_o}: S_{Q, z_o}\rightarrow \PP^1$ associated to $Q$ and $z_o$. Note that the covering transformation of $\widehat{S}_Q$ induces an involution on $S_{Q, z_o}$  which we will denote by $\sigma$.

 \[
\begin{CD}
\widehat{S}_{Q} @<{\mu}<< \overline{S}_{Q} @<{\nu_{z_o}}<<S_{\mcQ, z_o} \\
@V{\widehat{f}_{Q}}VV                 @VV{\overline{f}_{Q}}V      @ VV{\phi_{z_o}}V   \\
\PP^2@<<{q}< \overline{\PP^2}  @. \PP^1
\end{CD}
\]

We denote the set of sections of  $\phi_{z_o}$ by $\MW(S_{Q,z_o})$. The sections will be identified with their images and considered as curves on $S_{mcQ, z_o}$. It is known that  $\MW(S_{Q,z_o})$ can be endowed with an abelian group structure with a pairing $\langle, \rangle: \MW(S_{Q,z_o})\rightarrow Q$ called the {\it height pairing} (see \cite{shioda90}). When considering the height pairing, $ \MW(S_{Q,z_o})$ is called the Mordell-Weil lattice of $S_{Q,z_o}$.

Let $f=\widehat{f}_Q\circ\mu\circ\nu_{z_o}$. For a section $s\in\MW(S_{Q,z_o})$, let $C_s=f(s)$, the image of $s$ under $f$. The curve $C_s$  is a rational curve in $\PP^2$ whose local intersection numbers with $Q$ become even. Such curves are called contact curves of $Q$. Note that $f(s)=f(-s)$ where $-s$ is the negative of $s$ with respect to the group structure of $\MW(S_{Q,z_o})$. The curves $\mcC$ that we will study are reducible curves of the form
\[
\mcC=Q+C_{s_1}+\cdots+C_{s_r}
\]
for some choice of $s_1,\ldots, s_r\in \MW(S_{Q,z_o})$. The additional data related to  $\MW(S_{Q,z_o})$ allows us to distinguish the curves. 

Assume for simplicity that $\MW(S_{Q,z_o})$ is torsion free. Let  $E_i=\{s^i_1,\ldots, s^i_r\}\subset \MW(S_{Q,z_o})$ $(i=1,2)$ be subsets of $\MW(S_{Q,z_o})$ such that $C_{s_j^i}\not=C_{s_k^i}$ for $j\not=k$. We will consider the matroid structure on $E_1, E_2$ induced by the linear dependence relations in  $\MW(S_{Q,z_o})\otimes \mathbb{Q}$. Let $\mcC_i=Q+C_{s_1^i}+\cdots+C_{s_r^i}$ $(i=1,2)$.

\begin{thm}\label{thm:main1} 
Under the above settings, if $\MW(S_{Q,z_o})$ is torsion free and $E_1, E_2$ have distinct matroid structures, then there exist no homeomorphisms $h:\PP^2\rightarrow\PP^2$ with $h(\mcC_1)=\mcC_2$ and $h(Q)=Q$. 

Moreover, if $h(\mcC_1)=\mcC_2$ implies $h(Q)=Q$ necessarily and the combinatorics of $\mcC_1, \mcC_2$ are the same, then $(\PP^2,\mcC_1)$ and $ (\PP^2, \mcC_2)$  form a Zariski-pair.
\end{thm}


Theorem \ref{thm:main1}  allows us to distinguish Zariski pairs and Zariski $N$-ples by simply calculating the matroid structures of the subsets of $\MW(S_{Q,z_o})$. However, to actually construct Zariski pairs, we need to choose the subsets $\{s^i_1,\ldots, s^i_r\}$ so that they have the same combinatorics, which is a somewhat delicate matter. Fortunately, we were able to use classical results on smooth quartics and  bitangent lines, which can be found in \cite{dolg}, to overcome this difficulty. 

In the case where $Q$ is a smooth quartic, it is known that $\MW(S_{Q,z_o})\cong E_7^\ast$. The $E_7^\ast$ lattice has 28 pairs of minimal vectors  $\pm l_1,\ldots, \pm l_{28}$ of height $\frac{3}{2}$. Furthermore, $L_i=\mcC_{l_i}=\mcC_{-l_i}$ become bitangent lines of $Q$, and there is a bijection between the set of pairs $\pm l_i$ and the set of bitangent lines $L_i$. The combinatorics of these bitangent lines are known, as in the following proposition.

\begin{prop}\label{prop:combinat}
For a general smooth quartic $Q$, its bitangent lines $L_1,\ldots, L_{28}$ and a fixed value $r=1, \ldots, 28$, the combinatorics of curves of the form 
\[
Q+L_{i_1}+\cdots+L_{i_r}
\] 
are  the same for any $\{i_1,\ldots, i_r\}\subset\{1,\ldots,28\}$. Namely, all $L_{i_k}$ are true bitangents, i.e. they are tangent to $Q$ at two distinct points, and any three of $L_{i_1}, \ldots, L_{i_r}$ are non-concurrent.
\end{prop}

For curves $\mcC_1$, $\mcC_2$ of the form above, it is immediate that $h(\mcC_1)=\mcC_2$ implies $h(\mcQ)=\mcQ$ necessarily. Now,  Proposition \ref{prop:combinat} together with Theorem \ref{thm:main1} gives us the following theorem.

\begin{thm}\label{thm:main2}
Let $N_r$ be the number of distinct matroid structures on subsets of the form $\{l_{i_1},\ldots, l_{i_r}\}$, where $l_{i_k}$ is a representative of the pair $\pm l_{i_k}$.  Then there exists a Zariski $N_r$-ple of curves having the combinatorics as in Proposition \ref{prop:combinat}.
\end{thm}

At present, we have not been able to calculate the exact value of $N_r$ due to a lack of computer skills of the authors. However,  we have a lower bound as follows: 

\begin{prop}\label{prop:matroid}
For $r=1, \ldots, 28$, the value of $N_r$ is greater than or equal to $n_r$ given in the following table.
\begin{center}
\begin{tabular}{|c|c||c|c||c|c||c|c|}
\hline
$r$ & $n_r$ & $r$ & $n_r$ & $r$ & $n_r$ &  $r$ & $n_r$ \\
\hline
$1$ & $1$ & $8$ & $11$ & $15$ & $100$ & $22$ & $10$ \\
$2$ & $1$ & $9$ & $19$ & $16$ & $90$ & $23$ & $5$ \\
$3$ & $1$ & $10$ & $37$ & $17$ & $70$ & $24$ & $3$ \\
$4$ & $2$ & $11$ & $52$ & $18$ & $54$ & $25$ & $2$ \\
$5$ & $2$ & $12$ & $80$ & $19$ & $37$ & $26$ & $1$ \\
$6$ & $4$ & $13$ & $95$ & $20$ & $23$ & $27$ & $1$ \\
$7$ & $6$ & $14$ & $102$ & $21$ & $16$ & $28$ & $1$\\
\hline
\end{tabular}
\end{center}
\end{prop}

We remark that Zariski-pairs involving smooth quartics and its bitangent lines have already been studied by E. Artal-Bartolo and  J. Vall\`es.  They gave an example of a pair consisting of a smooth quartic and three bitangent lines. The results were privately communicated to the authors.  Also, the second author together with H. Tokunaga and M. Yamamoto have studied the case of  four bitangent lines where a Zariski triple exists. Our approach using matroids fails to detect these examples but we think that our work is still worthwhile as it is easy to increase the number of bitangent lines involved and can be applied to other non-smooth quartic curves. It also introduces a new point of view that is possibly relatively easier for a wider audience to access and hopefully will connect to other research areas.

The organization of this paper is as follows. In Section \ref{sec:prelim}, we review the basic terminology of matroids and results concerning elliptic surfaces and dihedral covers, which will give the connection between the matroid structure of sections and the topology of the curves. In Section 3, we will prove Theorem 1. In Section 4, we will discuss the  case where $Q$ is a smooth quartic and prove Theorem \ref{thm:main2} and also give the proof of Proposition \ref{prop:matroid}. In Appendix A, we give the source code used in our computations.

The second author is partially supported by Grant-in-Aid for Scientific Research C (18K03263).

\section{Preliminaries}\label{sec:prelim}
\subsection{Matroids}

As will be seen later, the (in)dependence of elements of $\MW(S_{Q,z_o})$ is deeply related to the (non)existence of certain Galois covers of $\PP^2$, hence it is important to understand the structure of (in)dependence. Here, Matroid Theory provides a nice framework as it was precisely designed to study generalizations of  the notion of linear independence in vector spaces.  In this section we briefly review the basic terminology of matroids. We refer to \cite{matroid} for more details.

There are many different cryptomorphic definitions of Matroids. In our paper, we are interested in the dependence of elements of $\MW(S_{Q,z_o})$, hence we adopt the definition based on {\it independent sets}. Let $E$ be a finite set and $2^E$ be the set of subsets of $E$.

\begin{defin}\label{def:matroid}
A matroid structure (or simply a matroid) on $E$ is a pair $(E, \mathcal{I})$, where $\mathcal{I}\subset 2^E$ satisfies
\begin{enumerate}
\item $\mathcal{I}\not=\emptyset$. (nontriviality)
\item For any $I_1, I_2\subset E$, if $I_1\subset I_2$ and $ I_2\in \mathcal{I}$, then $I_1\in\mathcal{I}$. (descending)
\item For every $I_1, I_2\in \mathcal{I}$, if $|I_1|<|I_2|$, then there exists $x\in I_2-I_1$ such that $I_1\cup \{ x\}\in \mathcal{I}$. (augmentation)  
\end{enumerate} 
Elements of $ \mathcal{I}$ will be called {\it independent sets} and the other subsets will be said to be {\it dependent}.
\end{defin}

\begin{example}
Let $V$ be a vector space, and $E=\{v_1, \ldots, v_r\}\subset V$. Let $\mathcal{I}=\{I\subset E \mid I \text{ is linearly independent} \}$. Then $\mathcal{I}$ clearly satisfies the conditions (1), (2), (3) in Definition \ref{def:matroid}. Hence $ (E, \mathcal{I})$ is a matroid structure on $E$. 
\end{example}

\begin{defin}
Let $(E,\mathcal{I})$ be a matroid. A subset $C\subset E$ is called a {\it circuit} if $C\not\in \mathcal{I}$ and all proper subsets of $C$ are independent sets. Moreover, $C$ is a minimal dependent set.
\end{defin}

\begin{example}
Let $V=\mathbb{R}^3$ and $v_1=\left(\begin{array}{c} 1 \\ 0 \\ 0 \end{array}\right)$, $v_2=\left(\begin{array}{c} 0 \\ 1 \\ 0 \end{array}\right)$, $v_3=\left(\begin{array}{c} 0 \\ 0 \\ 1 \end{array}\right)$ and $v_4=\left(\begin{array}{c} 1 \\ 1 \\ 1 \end{array}\right)$. Let $E=\{v_1, v_2, v_3, v_3\}$ and  consider the matroid structure induced by linear independence. Then $E$ itself forms a circuit.
\end{example}

%
%
%
%

\begin{defin}
Let $(E_1, \mathcal{I}_1)$, $(E_2, \mathcal{I}_2)$ be matroids. $(E_1, \mathcal{I}_1)$, $(E_2, \mathcal{I}_2)$ are said to be equivalent as matroids if there exists a bijection $\varphi:E_1\rightarrow E_2$ such that $I_1\in \mathcal{I}_1$ if and only if $\varphi(I_1)\in \mathcal{I}_2$. 
\end{defin}

\subsection{Elliptic surfaces and the Mordell-Weil lattice}

In this subsection, we list the basic facts about quartics, rational elliptic surfaces and the Mordell-Weil lattice.
We refer the reader to \cite{shioda90}, \cite{oguiso-shioda}  for more details. 

In this paper, an {\it elliptic surface} is a smooth projective surface $S$, with a relatively minimal genus 1 fibration $\phi:S\rightarrow C$ over a smooth projective curve $C$ having a section $O: C\rightarrow S$. We identify $O$  with its image in $S$. We also assume that $S$ has at least one singular fiber. Let ${\rm Sing}(\phi)=\{v\in C \mid \phi^{-1}(v) \text{ is singular }\}$. For $v\in {\rm Sing}(\phi)$, we put $F_v=\phi^{-1}(v)$ and denote its irreducible decomposition by $F_v=\Theta_{v,0}+\sum_{i=1}^{m_i-1} a_{v,i}\theta_{v,i}$, where $m_{v,i}$ is the number of irreducible components and $\Theta_{v,0}$ is the unique irreducible component with $\Theta_{v,0}.O=1$. The subset of ${\rm Sing}(\phi)$ that correspond to reducible singular fibers will be denoted by $R$. Let $\MW(S)$ be the set of sections of $\phi: S\rightarrow C$. 

The set $\MW(S)$ can be endowed with a group structure as follows. Let $E_S$ be the generic fiber of $\phi$ and $\mathbb{C} (C)$ be the function field of $C$. It is known that there is a bijection between  $\mathbb{C}(C)$ rational points $E_S(\mathbb{C}(C))$ of $E_S$  and $\MW(S)$. Furthermore, since we have $O \in\MW(S)$, $(E(S), O)$  can be considered as an elliptic curve over $\mathbb{C}(C)$ and has a group structure where $O$ acts as the identity element.

Furthermore, under these circumstances, $\MW(S)$ becomes a finitely generated abelian group with a pairing $ \langle,\rangle: \MW(S)\rightarrow \mathbb{Q}$ called the height pairing (\cite{shioda90}). 
The explicit formula to calculate the pairing for $s_1, s_2\in \MW(S)$  is given by
\[
\langle s_1, s_2 \rangle=\chi(S)+s_1.O+s_2.O-s_1.s_2-\sum_{v\in R} {\rm contr}_v(s_1,s_2) .
\] 
The formulas for calculating ${\rm contr}_v(s_1,s_2)$ can be found in \cite{shioda90}.

\subsection{Criterion for existence of dihedral covers}\label{sec:dihedral}

Let $D_{2n}$ be the dihedral group of order $2n$. We present a criterion for the existence of certain dihedral covers of $\PP^2$ in terms of $\MW(S)$. The existence/non-existence of the dihedral covers will enable us to distinguish the topology of the curves.

Let $Q$ be a quartic plane curve, $z_o\in Q$ be a general point of $Q$,  $s_1, \ldots, s_r\in\MW(S_{Q,z_o})$ be sections such that $C_{s_i}\not=C_{s_j}$, where $C_{s_i}=f(s_i)$ as in the Introduction.

\begin{thm}[{\cite[Corollary 4]{bannai-tokunaga15}}]\label{thm:dihedral}
Let $p$ be an odd prime. Under the above setting, there exists a $D_{2p}$-cover of $\PP^2$ branched at $2Q+p(C_{s_1}+\cdots+C_{s_{r}})$ if and only if there exists integers $a_i\in\{1,\ldots, p-1\}$ for $i=1, \ldots r$ such that $\sum_{ i=1} ^r a_i s_i\in p \MW(S)$.
\end{thm}

\begin{cor}\label{cor:dependence}
If there exists a $D_{2p}$ cover branched at $2Q+p(C_{s_1}+\cdots+C_{s_{r}})$, then the images of $s_1, \ldots, s_r$ in $\MW(S)\otimes \ZZ/p\ZZ$ become linearly dependent.
\end{cor}

Note that the converse of Corollary \ref{cor:dependence}  is not true, as it is necessary for the images of $s_1, \ldots, s_r$ to have a linear dependence relation where all coefficients are non-zero for there to exist a dihedral cover. If there does not exist such linear dependence relation, the branch locus will not be the whole of $2Q+p(C_{s_1}+\cdots+C_{s_{r}})$. To exclude such cases, the notion of  circuits is useful.

\begin{cor}\label{cor:circuit}
If the images of $s_1,\ldots, s_r$ in $\MW(S)\otimes \ZZ/p\ZZ$ forms a circuit, then there exists a $D_{2p}$-cover branched at $2Q+p(C_{s_1}+\cdots+C_{s_{r}})$.
\end{cor}

If $s_1,\ldots, s_r$ form a circuit over $\QQ$, then their images in $\MW(S)\otimes\ZZ/p\ZZ$  form a circuit for infinitely many prime numbers $p$. Hence we have:

\begin{lem}\label{cor:existance}
If $s_1,\ldots, s_r$ are linearly dependent, then there are infinitely many prime numbers $p$ such that there exists a $D_{2p}$-cover branched at $2Q+p(C_{s_{i_1}}+\cdots+C_{s_{i_t}})$ for some nonempty subset $\{i_1, \ldots, i_t\}\subset \{1,\ldots, r\}$. 

\end{lem}

On the other hand, if $s_1,\ldots, s_r$ are independent over $\QQ$, then they are independent over $\ZZ/p\ZZ$ except for a finite number of primes. This implies the following.

\begin{lem}\label{lem:independence}
If $s_1, \ldots, s_r$ are independent over $\QQ$, then there are only a finite number of prime numbers $p$ such that there exists a $D_{2p}$-cover branched at $2Q+p(C_{s_{i_1}}+\cdots+C_{s_{i_t}})$ for some nonempty subset $\{i_1, \ldots, i_t\}\subset \{1,\ldots, r\}$. 
\end{lem}

\section{Proof of Theorem 1}\label{sec:proof1}
In this section, we use the criterion for the existence of dihedral covers given in Section \ref{sec:dihedral} to connect the data of matroids of subsets of $\MW(S_{\mcQ, z_o})$ to the  data of the embedded topology of the curves in $\PP^2$, and prove Theorem \ref{thm:main1}.

Let $E_i=\{s_1^i,\ldots, s_r^i\}\subset \MW(S_{Q, z_o})$ $(i=1,2)$ be subsets of $\MW(S_{\mcQ, z_o})$ such that $C_{s_j^i}\not =C_{s_k^i}$ for $j\not=k$. Consider the matroid structure $(E_i, \mathcal{I}_i)$ on $E_i$ $(i=1,2)$ induced by the linear dependence relation in $\MW(S_{Q,z_o})\otimes \mathbb{Q}$. Let $\mcC_i=Q+C_{s_1^i}+\cdots+C_{s_r^i}$ $(i=1,2)$.

\begin{prop}\label{prop:key}
If there exists a homeomorphism $h:\PP^2\rightarrow \PP^2$ such that $h(\mcC_1)=\mcC_2$ and $h(Q)=Q$, then $(E_1, \mathcal{I}_1)$ and $(E_2, \mathcal{I}_2)$ are equivalent as matroids.
\end{prop}
\proof
By the assumption that $h$ is a homeomorphism such that $h(\mcC_1)=\mcC_2$ and $ h(Q)=Q$, $h$ induces a bijection $\{C_{s_1^1},\ldots, C_{s_r^1}\}\rightarrow  \{C_{s_1^2},\ldots, C_{s_r^2}\}$ which in turn induces a bijection $h_\ast: E_1\rightarrow E_2$. Let $I_1\in\mathcal{I}_1$ be an independent set. Then by  Lemma \ref{lem:independence}, there exists only a finite number of primes such that a  $D_{2p}$ cover branched at $2Q+p(\sum_{s\in J_1} C_{s})$ for some subset $J_1\subset I_1$ exists.  Since $h$ is a homeomorphism, the same is true for $h_\ast(I_1)$ which implies that $h_\ast(I_1)\in\mathcal{I}_2$,  by Lemma \ref{cor:existance}. The converse is also true so we have $I_1\in \mathcal{I}_1$ if and only if $h_\ast(I_1) \in \mathcal{I}_2$.  Therefore $(E_1, \mathcal{I}_1)$ and $(E_2, \mathcal{I}_2)$ are equivalent as matroids.
\qed

The contrapositive of Proposition \ref{prop:key} gives Theorem \ref{thm:main1}. 

\begin{remark} 
The statement of Proposition \ref{prop:key} concerns the matroid structure over $\mathbb{Q}$. However, from the proof,  it is evident that if we consider the matroid structures of the sections in $\MW(S) \otimes \ZZ/p\ZZ$ for all $p$ we would be able to distinguish the arrangements in more detail. 
\end{remark}

\section{The smooth case}\label{sec:smoothcase}

In this section, we will consider the case where $Q$ is a smooth quartic.

\subsection{The bitangents of $Q$ and sections of $S_{Q,z_o}$}

We will use the notation given in the Introduction. Let $Q$ be a smooth plane quartic and $z_o\in Q$ be a general point of $Q$. Since $Q$  is smooth, $\widehat{S}_{Q}=\overline{S}_Q$. In this case $S_{Q,z_o}$ has only one reducible  singular fiber $F_0=\Theta_{0,0}+\Theta_{0,1}$ of type $\I_2$. The component $\Theta_{0,0}$ is the exceptional divisor of the first blow up of $\mu_{z_o}$  in the introduction, and $\Theta_{0,1}$ is the strict transform of the preimage of the tangent line of $Q$ at $z_o$. All other singular fibers are irreducible. By \cite{oguiso-shioda}, we have $\MW(S_{Q,z_o})\cong E_7^\ast$ where $E_7^\ast$ is the dual lattice of the root lattice $E_7$. 
It is known that the $E_7^\ast$ lattice has 56 minimal vectors $\pm {l_1},\ldots, \pm l_{28}$ of height $\frac{3}{2}$. It is also well known that $Q$  has 28 bitangent lines $L_1,\ldots, L_{28}$. The correspondence between the 28 pairs of minimal vectors and  the 28 bitangent lines is given in \cite{shioda93}, but we describe the relation below for the readers convenience. 

\begin{lem}\label{lem:vec2bitan}
Let $l\in \MW(S_{Q,z_o})$ be a minimal vector of height $\frac{3}{2}$. Then $L=f(l)$ is a bitangent line of $Q$, where $f$ is the morphism $f: S_{Q,z_o}\rightarrow \PP^2$ given in the Introduction.
\end{lem}  
\proof
By the explicit formula for the height pairing, and since $\chi(S_{Q,z_o})=1$ and $l.l=-1$, we have
\[
\langle l, l\rangle=2+2l.O-{\rm contr}(l,l)=\frac{3}{2}.
\]
Where ${\rm contr}(l,l)$ is the contribution from the unique reducible singular fiber $F_0$. Since the possible values of ${\rm contr}(l,l)=0, \frac{1}{2}$, we have $l.O=0$ and ${\rm contr}(l,l)=\frac{1}{2}$ which implies that $l.\Theta_{0,1}=1$. This implies that $l$ is disjoint with the exceptional set of $\nu_{z_o}$. 
Also, if we consider the section $-l=\sigma^\ast(l)$, the preimage of $l$ under the involution $\sigma$, we have
\[
\langle l, -l\rangle=1+l.O+(-l).O-l.(-l)-{\rm contr}(l,-l)=-\frac{3}{2}
\]
Hence we obtain $l.(-l)=2$. Let $\widehat{l}=\nu_{z_o}(l)$ and $\widehat{-l}=\nu_{z_o}(-l)$. The above implies that $\widehat{l}.\widehat{-l}=\widehat{l}.\widehat{Q}=2$, where $\widehat{Q}$ is the ramification locus of $\widehat{f}_Q$. Now since $(\widehat{f}_Q)^\ast(L)=\widehat{l}+\widehat{-l}$ we have  $2L.L=(\widehat{l}+\widehat{-l}).(\widehat{l}+\widehat{-l})$. Hence we obtain $L.L=1$ which implies that $L$ is a line in $\PP^2$. Also, the local intersection numbers of $L$ and $Q$ must be even by construction, hence $L$ is a bitangent line. 
\qed

\begin{rem}
Note that the two points of tangency may coincide to give a line $L$  intersecting $Q$ at a single point with multiplicity 4, which we will still consider to be a bitangent line.
\end{rem}
\begin{lem}
Let $L$ be a bitangent line of $Q$ and let $f^\ast(L)=l+l^\prime$. Then $l, l^\prime$ become minimal sections with height $\frac{3}{2}$ and $l^\prime =\sigma^\ast{l}=-l$.
\end{lem}
\proof 
By following through the proof of Lemma \ref{lem:vec2bitan} backwards, we have the desired result.
\qed

The above two lemmas give us the following propositon.
\begin{prop}
There is a bijection between the set of  28 bitangent lines of $Q$ and the set of 28 pairs of minimal vectors of the $E_7^\ast$ lattice.
\end{prop}

\subsection{Riemann's Equations for bitangents}

In this subsection we prove Proposition \ref{prop:combinat} by using Riemann's Equations for bitangents of $Q$. The details about Riemann's Equations including the proofs and historical notes can be found in \cite{dolg}. However, the equations given there have some typos so we will restate the correct equations here for the readers convenience.

Given the equation of seven btangent lines $L_1,\ldots, L_7$ of $Q$, which form an Aronhold set, it is possible to recover the defining equation of $Q$ and the equations of the remaining 21 biitangent lines. We can assume that $L_1, \ldots, L_7$ are given by the following equations for a suitable choice of coordinates:

\begin{eqnarray*}
L_1=V(t_0), L_2=V(t_1), L_3=V(t_2), L_4=V(t_0+t_1+t_2) \\
L_{4+i}=V(a_{0i}t_0+a_{1i}t_1+a_{2i}t_2), (i=1,2,3)
\end{eqnarray*}

\begin{thm}[\cite{dolg}, Theorem 6.1.9]
There exists linear forms $u_0, u_1, u_2$ such that, after rescaling the forms, 
\[
C=V(\sqrt{t_0u_0}+\sqrt{t_1u_1}+\sqrt{t_2u_2}).
\]
The forms $u_0, u_1, u_2$ can be found from equations
\begin{eqnarray*}
u_0+u_1+u_2+t_0+t_1+t_2=0\\
\frac{u_0}{a_{01}}+\dfrac{u_1}{a_{11}}+\frac{u_2}{a_{21}}+k_1a_{01}t_0+k_1a_{11}t_1+k_1a_{21}t_2=0\\
\frac{u_0}{a_{02}}+\dfrac{u_1}{a_{12}}+\frac{u_2}{a_{22}}+k_2a_{02}t_0+k_2a_{12}t_1+k_2a_{22}t_2=0\\
\frac{u_0}{a_{03}}+\dfrac{u_1}{a_{13}}+\frac{u_2}{a_{23}}+k_3a_{03}t_0+k_3a_{13}t_1+k_3a_{23}t_2=0
\end{eqnarray*}
where $k_1, k_2, k_3$ can be found from solving first, 
\[
\left(\begin{array}{ccc} 
\frac{1}{a_{01}} & \frac{1}{a_{02}} & \frac{1}{a_{03}} \\
\frac{1}{a_{11}} & \frac{1}{a_{12}} & \frac{1}{a_{13}} \\
\frac{1}{a_{21}} & \frac{1}{a_{22}} & \frac{1}{a_{23}} 
 \end{array}\right)\left(\begin{array}{c} \lambda_1 \\ \lambda_2 \\ \lambda_3 \end{array}\right)= \left(\begin{array}{c} -1 \\ -1 \\ -1 \end{array}\right)
\]
and then solving
\[
\left(\begin{array}{ccc} 
\lambda_1a_{01} & \lambda_2a_{02} & \lambda_3 a_{03} \\
\lambda_1a_{11} & \lambda_2a_{12} & \lambda_3 a_{13} \\
\lambda_1a_{21} & \lambda_2a_{22} & \lambda_3 a_{23} 
 \end{array}\right)\left(\begin{array}{c} k_1 \\ k_2 \\ k_3 \end{array}\right)= \left(\begin{array}{c} -1 \\ -1 \\ -1 \end{array}\right).
\]
\end{thm}

 \begin{thm}[\cite{dolg}, Theorem 6.1.9]
Given the equations of an Aronhold set as in the previous Theorem, the equations of the remaining 21 bitangent lines are given by
\begin{enumerate}
\item $u_0=0$, $u_1=0$, $u_2=0$
\item $u_0+t_1+t_2=0$, $t_0+u_1+t_2=0$, $t_0+t_1+u_2=0$
\item $\dfrac{u_0}{a_{0i}}+k_i(a_{1i}t_1+a_{2i}t_2)$, $(i=1,2,3)$
\item $\dfrac{u_1}{a_{1i}}+k_i(a_{0i}t_0+a_{2i}t_2)$, $(i=1,2,3)$
\item $\dfrac{u_2}{a_{2i}}+k_i(a_{0i}t_0+a_{1i}t_1)$, $(i=1,2,3)$
\item $\dfrac{t_0}{1-k_ia_{1i}a_{2i}}+\dfrac{t_1}{1-k_ia_{0i}a_{2i}}+\dfrac{t_2}{1-k_ia_{0i}a_{1i}}$, $(i=1,2,3)$
\item $\dfrac{u_0}{a_{0i}(1-k_ia_{1i}a_{2i})}+\dfrac{u_1}{a_{1i}(1-k_ia_{0i}a_{2i})}+\dfrac{u_2}{a_{2i}(1-k_ia_{0i}a_{1i})}$, $(i=1,2,3)$
\end{enumerate} 
\end{thm}

Since we have explicit equations, it is possible to calculate the combinatorics of the bitangent lines. We used the open-source mathematics software system  SageMath \cite{sage} for the actual calculations.

\begin{lem}\label{lem:combinat}

For a general choice of $a_{0i}, a_{1i}, a_{2i}$ $(i=1,2,3)$ the following hold:
\begin{enumerate}
\item Any three bitangent lines of $Q$ are non-concurrent.
\item Every bitangent line of $Q$ is a true bitangent, i.e. it is tangent to $Q$ at two distinct points.
\end{enumerate}
\end{lem}
\proof
Since the condition for three lines to be concurrent is a closed condition on $a_{0i}, a_{1i}, a_{2i}$ $(i=1,2,3)$, it is enough to find one example where the statement holds. Almost any choice will serve our purpose. The same is true for the second statement. 
\qed

Lemma \ref{lem:combinat} immediately implies Proposition \ref{prop:combinat}.

\subsection{The proof of Proposition \ref{prop:matroid}}

In this subsection, we describe the method we used to distinguish the matroid structures of minimal vectors of the $E_7^\ast$ lattice in order to calculate $n_r$. We used SageMath \cite{sage} for the actual calculations. 

The object that we want to classify are the matroid structures on the sets of the form $\{l_{i_1},\ldots, l_{i_r}\}$ where $l_{i_r}$ are representatives of pairs $\pm l_{i_r}$ of minimal vectors of height $\frac{3}{2}$.
It is known that the $E_7^\ast$ lattice can be representation in $\mathbb{Q}^8$ in a way so that the minimal vectors are of the form
\[
\pm\frac{1}{4}(1,1,1,1,1,1,-3,-3)
\]  
and its permutations. We use this representation in our calculations.

We used an inductive argument on the number of vectors $r$. 
For each subset $E\subset \{l_1,\ldots, l_{28}\}$ having $r$-elements, we assign an $(n_{r-1}+1)$-ple of integers inductively as follows. The values of $n_{r}$ will also be determined inductively along the way.
\begin{itemize}
\item {\bf Step $(1)$}

For every subset with a single element, we assign the pair $\alpha_{1,1}=(1;1)$.

\item {\bf Step $(k+1)$}

Suppose that every subset having $k$ elements has been assigned an $(n_{k-1}+1)$-ple of integers.  We set $n_k$ to be the number of distinct $(n_{k-1}+1)$-ples that have been assigned and label them by $\alpha_{k,1},\ldots, \alpha_{k,n_k}$. Next, to each subset $E\subset \{l_1,\ldots, l_{28}\}$  having $k+1$ elements, we assign an $(n_k+1)$-ple as follows:
\begin{enumerate}
\item[(i)] Consider the linear dependence/independence of $E$. Put $i=0$ if it is dependent and $i=1$ if it is independent.
\item[(ii)] Let $m_j^k$ be the number of subsets of $E$ of $k$ elements that have the  $(n_{k-1}+1)$-ple $\alpha_{k,j}$  assigned.
\item[(iii)] Assign the $(n_k+1)$-ple $ (i; m_1^k, \ldots, m_{n_k}^k)$ to $E$.
\end{enumerate}
\end{itemize}

\begin{lem}\label{lem:lbd}
Let $E_1$, $E_2$ be subsets of $\{ l_1, \ldots, l_{28}\}$  and $|E_1|=|E_2|=r$. If $E_1$ and $E_2$ have the same matroid structure, then the $(n_{r-1}+1)$-ples of integers assigned above are equivalent. 
\end{lem}
\proof  We use induction on $r$ to prove this lemma. The case for $r=1$ is trivial as every subset having a single element has the same pair assigned and has the same matroid structure.

Assume the statement holds for $r=k$. If $|E_1|=|E_2|=k+1$ and $E_1, E_2$ have equivalent matroid structure, there exists a bijection $\varphi: E_1\rightarrow E_2$ that preserves independent sets. Hence $E_1$ is independent if and only if $E_2$ is independent and the value of $i$ must be equal. Also, $\varphi$ induces a bijection from $\{E\subset E_1 \mid |E|=k\}$ to $\{E\subset E_2 \mid |E|=k\}$ and an equivalence of matroid structures among the corresponding subsets. Hence the values of $m_j^k$ must be equal do to the hypothesis of induction, and the assigned $(n_k+1)$-ple are equivalent.
 \qed
 
 Lemma \ref{lem:lbd} and calculations done by computer using SageMath gives Proposition \ref{prop:matroid}.

\noindent Ryutaro SATO, Shinzo BANNAI\\
National Institute of Technology, Ibaraki College\\
866 Nakane, Hitachinaka-shi, Ibaraki-Ken 312-8508 JAPAN \\
{\tt sbannai@ge.ibaraki-ct.ac.jp}\\
\newpage

\appendix
\section{Implementation}

In this appendix, we show the source code we implemented for distinguishment. In practice, we devised to reduce the amount of space complexity to execute this program.

\begin{lstlisting}
from itertools import combinations
from itertools import permutations

def combine(a):
    A = []
    combine2(A, a, 0)
    return A

def combine2(A, a, n):
    if(len(a) == n):
        A.append(list(a))
        return

    a[n] = 1
    combine2(A, a, n + 1)

    a[n] = -3
    combine2(A, a, n + 1)

L = [l for l in sorted(set(list(permutations([1, 1, 1, 1, 1, 1, -3, -3]))), reverse = True)]
row = {tuple(l) : i + 1 for i, l in enumerate(L)}

ind = [{} for i in range(29)]
dep = [{} for i in range(29)]
n_r = [{} for i in range(29)]

for n in range(1, 29):
    if n < 8:
        col = {tuple(l) : i for i, l in enumerate(combine([-3 for i in range(n)]))}
        col_inv = {i : tuple(l) for l, i in col.items()}

        row_com = [l for l in combinations(L, r = n)]
        col_com = [list(col[tuple([l[i][j] for i in range(n)])] for j in range(8)) for l in combinations(L, r = n)]

        uni = [list(y) for y in set(tuple(sorted(x)) for x in col_com)]
        mat = [matrix(QQ, [col_inv[l[i]] for i in range(8)]) for l in uni]
        res = set([tuple(uni[i]) for i, j in enumerate(mat) if rank(j) == n])

        ind[n] = set([tuple([row[tuple(x[j])] for j in range(n)]) for x in ([row_com[i] for i, y in enumerate(col_com) if(tuple(sorted(y)) in res)])])
        dep[n] = set(x for x in set(tuple(row[tuple(y)] for y in z) for z in row_com) if x not in ind[n])

    else:
        row_i = [i for i in range(1, 29)]
        com = [tuple(i) for i in combinations(row_i, r = n)]

        ind[n] = {}
        dep[n] = set(tuple(com))

    if n == 1:
        continue

    if len(dep[n]) == 0:
        n_r[n] = {x:tuple([0]) for x in ind[n]}

    elif len(ind[n]) != 0:
        n_r[n] = {x:tuple(sorted([n_r[n - 1][y] for y in combinations(x, r = n - 1) if y in n_r[n - 1]])) for x in dep[n]}
        for x in ind[n]:
            n_r[n][x] = tuple([0])

    else:
        n_r[n] = {x:tuple(sorted([n_r[n - 1][y] for y in combinations(x, r = n - 1)])) for x in dep[n]}

    print("r={0},n_r={1}".format(n, len(set(n_r[n].values()))))


\end{lstlisting}

\end{document}